\documentclass[10pt, reqno]{amsart}
\usepackage{amsmath, amsthm, amscd, amsfonts, amssymb, graphicx, color}
\usepackage[bookmarksnumbered, colorlinks, plainpages]{hyperref}
\hypersetup{colorlinks=true,linkcolor=red, anchorcolor=green, citecolor=cyan, urlcolor=red, filecolor=magenta, pdftoolbar=true}
\usepackage{mathrsfs}

\newtheorem{theorem}{Theorem}[section]

\newtheorem{corollary}[theorem]{Corollary}
\theoremstyle{definition}
\newtheorem{definition}[theorem]{Definition}
\newtheorem{example}[theorem]{Example}

\theoremstyle{remark}
\newtheorem{remark}[theorem]{Remark}
\numberwithin{equation}{section}

\begin{document}
\setcounter{page}{1}

\title[Noncommutative geometry of elliptic surfaces]{Noncommutative geometry of elliptic surfaces}

\author[Nikolaev]
{Igor V. Nikolaev$^1$}

\address{$^{1}$ Department of Mathematics and Computer Science, St.~John's University, 8000 Utopia Parkway,  
New York,  NY 11439, United States.}
\email{\textcolor[rgb]{0.00,0.00,0.84}{igor.v.nikolaev@gmail.com}}


\subjclass[2010]{Primary 11G05; Secondary 46L85.}

\keywords{elliptic surfaces, continued fractions.}


\begin{abstract}
We recast elliptic surfaces over the projective line in terms of the 
non-commutative tori and one-parameter families of the periodic continued fractions.   
The correspondence is used to study the Picard numbers, the ranks and  the minimal models  of 
such surfaces. As an example, we calculate the Picard numbers 
of elliptic surfaces with fibers having complex multiplication.  
\end{abstract}

\maketitle

\section{Introduction}
Our paper deals with two types of objects: the well known -- elliptic curves 
and surfaces, and the less known -- non-commutative tori and certain families
of periodic continued fractions.  By magic and glory of mathematics,  they 
happen to be the same, the second (less known) exposing an intrinsic structure, 
while the first (well known) following a historic pattern. 
We hope that the reader will find the second type  more natural 
and accurate description of the elliptic curves and surfaces.

Recall that elliptic surfaces are bundles over a  base curve of genus $g\ge 0$
with   fibers of  genus one. 
We assume further that $g=0$  and denote by  $\mathbf{Q}(t)$  the 
field of rational functions on the base $\mathbf{Q}P^1$. 
The elliptic curve over the field  $\mathbf{Q}(t)$ can be identified with an 
elliptic surface [Sch\"utt \& Shioda 2019] \cite[Chapter 5]{SS} which we
write $\mathscr{E}_{\mathbf{Q}(t)}$ by an abuse of notation.   
The  generic fiber of  $\mathscr{E}_{\mathbf{Q}(t)}$
will be denoted by $\mathscr{E}_{t}$. 
Silverman's specialization theorem says that the  homomorphism 
$\mathscr{E}_{\mathbf{Q}(t)}\to \mathscr{E}_{t}$ is injective for all but finitely many rational values of $t$. 
In particular,  the rank  of $\mathscr{E}_{t}$ is related to the rank of 
$\mathscr{E}_{\mathbf{Q}(t)}$.

The less known (but equivalent) definition of the fiber  $\mathscr{E}_{t}$
is as follows. The Sklyanin algebra $S(\alpha,\beta,\gamma)$ is a free $\mathbf{C}$-algebra on  four generators
$\{x_1, \dots, x_4\}$ satisfying six quadratic relations
$x_1x_2-x_2x_1 = \alpha(x_3x_4+x_4x_3)$,
\quad $x_1x_2+x_2x_1 = x_3x_4-x_4x_3$,
~$x_1x_3-x_3x_1 = \beta(x_4x_2+x_2x_4)$,
~$x_1x_3+x_3x_1 = x_4x_2-x_2x_4$,
~$x_1x_4-x_4x_1 = \gamma(x_2x_3+x_3x_2)$ and 
$x_1x_4+x_4x_1 = x_2x_3-x_3x_2$,
where $\alpha,\beta,\gamma\in \mathbf{C}$ and  
$\alpha+\beta+\gamma+\alpha\beta\gamma=0$. 
The algebra $S(\alpha,\beta,\gamma)$ is  twisted homogeneous coordinate ring  of an elliptic curve in the Jacobi form  $\mathscr{E}\subset \mathbf{C}P^3$,
i.e.  as an intersection of the quadrics $u^2+v^2+w^2+z^2 ={1-\alpha\over 1+\beta}v^2+{1+\alpha\over 1-\gamma}w^2+z^2 = 0$  
 [Stafford \& van ~den ~Bergh 2001]  \cite[Example 8.5]{StaVdb1}.
The norm-closure of a self-adjoint representation of the algebra $S(\alpha,\beta,\gamma)$ by 
bounded linear operators on a Hilbert   space $\mathscr{H}$  is a  non-commutative torus $\mathscr{A}_{\theta}$,   
 i.e.   a $C^*$-algebra  generated by a pair of unitary operators $u$ and $v$
satisfying the commutation relation $vu=e^{2\pi i\theta}uv$  for a real constant $\theta$.
The $\mathscr{A}_{\theta}$ is said to have real multiplication (RM), if $\theta$ 
is a quadratic irrationality given by   $k$-periodic continued 
fraction $[b_1,\dots, b_N;  \overline{a_1,\dots,a_k}]$. 
The map $F: \mathscr{E}_t\mapsto \mathscr{A}_{\theta}$ is a functor, 
such that if $t\in\mathbf{Q}$,  then $F(\mathscr{E}_t)$ are  non-commutative
 tori with RM  \cite[Section 1.3]{N}. 
\begin{example}\label{exm1.2}
{\bf (\cite{Nik1})}
Let  $\mathscr{E}_{\mathbf{Q}(t)}$ be a surface
 given by the affine equation:
\begin{equation}\label{eq1.5}
y^2=x(x-1)\left(x-{t-2\over t+2}\right),\quad t\in\{3, 4, 5, \dots\}.
\end{equation}
Then  $F$ acts by the formula: 
\begin{equation}\label{eq1.6}
\mathscr{E}_{\mathbf{Q}(t)}\mapsto \mathscr{A}_{[t-1;  ~\overline{1, ~t-2}]}.
\end{equation}
\end{example}

\bigskip
The aim of our note is a generalization of  (\ref{eq1.6}) to the rest of elliptic surfaces $\mathscr{E}_{\mathbf{Q}(t)}$
 in terms of  the continued fractions $[b_1,\dots, b_N; \overline{a_1,\dots,a_k}]$,
 see Theorem \ref{thm1.1}.  The result is used   to calculate
the Picard number, the rank and  the minimal model for the $\mathscr{E}_{\mathbf{Q}(t)}$,
see Theorem \ref{thm1.3}. The following notation will be used.

\begin{definition}
{\bf (Brock, Elkies \& Jordan \cite{BrElJo1})}
We denote by $V_{N,k}(\mathbf{C})$ an affine variety defined by 
the polynomials in variables $a_{1\le i\le k}$ and $b_{1\le j\le N}$
satisfying the following obvious equality: 
\begin{equation}\label{eq1.1}
[b_1,\dots, b_N, a_1,\dots, a_k; \overline{a_1,\dots,a_k}]=[b_1,\dots, b_N; \overline{a_1,\dots,a_k}].
\end{equation}
 \end{definition}
\begin{remark}
The equations of the
 $V_{N,k}(\mathbf{C})$ date back to  [Euler 1765] \cite{Eul1}. 
The corresponding projective variety was studied  in \cite[Section 6.2.1]{N}. 
Our notation $V_{N,k}(\mathbf{C})$ corresponds to the variety $V(\mathcal{B})_{N,k}$ 
of [Brock, Elkies \& Jordan 2021]  \cite{BrElJo1}. 
\end{remark}
The integer points of  $V_{N,k}(\mathbf{C})$ parametrize the $k$-periodic  continued 
fractions $[b_1,\dots, b_N; \overline{a_1,\dots,a_k}]$
and  $\dim ~V_{N,k} (\mathbf{C})=N+k-2$.  The  $V_{N,k}(\mathbf{C})$ 
is a fiber bundle over the Fermat-Pell conic  $\mathscr{P}:~Cx^2-Bxy+Ay^2=(-1)^kA$
with the fiber map $\pi: V_{N,k}(\mathbf{C})\to\mathscr{P}$, 
see  [Brock, Elkies \& Jordan 2021] \cite{BrElJo1}
or (\ref{eq2.8}) for the details.  
By $t\in\mathbf{C}P^1$ we denote a rational parametrization of $\mathscr{P}$.
In what follows, we consider only regular maps, i.e. given by polynomials over
the base field. 
Our main result can be formulated as follows.
\begin{theorem}\label{thm1.1}
For each elliptic surface $\mathscr{E}_{\mathbf{Q}(t)}$
there exists  a section  
$\mathscr{P}\to$ 
\linebreak
$U_{b_1,\dots,b_N; a_1,\dots, a_k}$ of  a  sub-bundle 
$(U_{b_1,\dots,b_N; a_1,\dots, a_k}, \mathscr{P}, \pi')$ of the fiber bundle 
\linebreak
$(V_{N,k}(\mathbf{C}), \mathscr{P}, \pi)$, 
 such that:

\begin{equation}\label{eq1.4}
F\left(\mathscr{E}_{\mathbf{Q}(t)}\right)=\left\{\mathscr{A}_{[b_1(t),\dots, b_N(t);  ~\overline{a_1(t),\dots,a_k(t)}]}
~|~a_i(t), b_j(t)\in  \mathbf{N},
~t\in\mathscr{P}
\right\}.
\end{equation}
\end{theorem}

  \bigskip
  Denote by $\mathscr{F}=(\pi')^{-1}(t)$ a fiber of the bundle $(U_{a_1,\dots,a_N; b_1,\dots, b_k}, \mathscr{P}, \pi')$
  over $t\in \mathscr{P}$. 
 For  $p(t)\in\mathbf{Z}[t]$
 let  $\mathscr{E}^{CM}_{p(t)}$ be  an elliptic surface 
  whose fibers over $t\in\mathbf{Q}$  have  complex multiplication (CM) by  
 the ring of integers of the imaginary quadratic field  $\mathbf{Q}\left(\sqrt{-1-p^2(t)}\right)$,
 see Section 4 for an example.  
Denote by  $\mathscr{E}^{\min}_{\mathbf{Q}(t)}$ the minimal model of 
the surface $\mathscr{E}_{\mathbf{Q}(t)}$.    An application of  theorem \ref{thm1.1} is as follows. 
\begin{theorem}\label{thm1.3}
 Let  $\mathscr{E}_{\mathbf{Q}(t)}$ be a surface  satisfying  equation  (\ref{eq1.4}) and   $\rho(\mathscr{E}_{\mathbf{Q}(t)})$
 its Picard number over $\mathbf{Q}$.  The following is true:
 
 \medskip
 (i) the Picard number $\rho(\mathscr{E}_{\mathbf{Q}(t)})=N+k$; 

\smallskip
(ii) the rank $r(\mathscr{E}_{\mathbf{Q}(t)})=\dim ~\mathscr{F}$;

\smallskip
(iii)  the $\mathscr{E}^{\min}_{\mathbf{Q}(t)}$ is 
$\mathbf{C}$-isomorphic to the $\mathscr{E}^{CM}_{p(t)}$ 
for some $p(t)>0$. 
\end{theorem}
The article is organized as follows. Preliminary facts can be found 
in Section 2.  The proofs of theorems \ref{thm1.1} and \ref{thm1.3} 
 are given in Section 3.  The Picard numbers of
 the elliptic surfaces with generic fiber having complex multiplication
 are calculated in Section 4.

\section{Preliminaries}
We briefly review continued fractions,  elliptic surfaces  and non-commutative tori. 
For a detailed exposition we refer the reader to [Brock, Elkies \& Jordan 2021] \cite{BrElJo1},
[Sch\"utt \& Shioda 2019] \cite[Chapter 5]{SS} and \cite[Chapter 1]{N}, respectively.

\subsection{Continued fractions}
By an infinite continued fraction one understands an expression of the form:
\begin{equation}\label{eq2.1}
[c_1,c_2, c_3, \dots]:=c_1+\cfrac{1}{c_2+\cfrac{1}{c_3+\dots}} ~,
\end{equation}
where $c_1$ is an integer and $c_2,c_3,\dots$ are positive integers. 
The continued fraction (\ref{eq2.1}) converges to an irrational number and each irrational 
number has a unique representation by (\ref{eq2.1}).  The expression (\ref{eq2.1})
is called $k$-periodic, if $c_{i+k}=c_i$ for all $i\ge N$ and a minimal index $k\ge 1$.   
We shall denote the $k$-periodic continued fraction by 
\begin{equation}\label{eq2.2}
[b_1,\dots, b_N, \overline{a_1,\dots,a_k}],
\end{equation}
 where $(a_1,\dots, a_k)$ is the minimal period of (\ref{eq2.1}). 
The continued fraction (\ref{eq2.2}) converges to one of the irrational root of a quadratic 
polynomial 
\begin{equation}\label{eq2.3}
Ax^2+Bx+C\in \mathbf{Z}[x].
\end{equation}
Conversely,  the irrational root of any quadratic polynomial  (\ref{eq2.3}) has a
representation by  the continued fraction (\ref{eq2.2}). 
Notice that the following two continued fractions define the same irrational number: 
\begin{equation}\label{eq2.4}
[b_1,\dots, b_N, \overline{a_1,\dots,a_k}]= [b_1,\dots, b_N, a_1,\dots, a_k, \overline{a_1,\dots,a_k}]. 
\end{equation}
But it is well known,  that two infinite continued fraction with at most finite number of distinct 
entries must be related by the linear fractional transformation given  by a matrix  $\mathcal{E}\in GL_2(\mathbf{Z})$.
Therefore equation (\ref{eq2.4}) can be written in the form
\begin{equation}\label{eq2.5}
x={E_{11}x+E_{12}\over E_{21}x+E_{22}},
\end{equation}
where   $\mathcal{E}=(E_{ij})\in GL_2(\mathbf{Z})$ and $x=[b_1,\dots, b_N, \overline{a_1,\dots,a_k}]$. 
\begin{remark}\label{rmk2.1}
It is easy to see, that $x$ in  (\ref{eq2.5}) is  the root of  quadratic polynomial (\ref{eq2.3}) with $A=E_{21}, B=E_{22}-E_{11}$
and $C=-E_{12}$. 
\end{remark}
\begin{definition}\label{dfn2.2}
The  Brock-Elkies-Jordan variety  $V_{N,k}(\mathbf{C})\subset \mathbb{A}^{N+k}$ is an affine variety over $\mathbf{Z}$  
defined by the three equations:
\begin{equation}\label{eq2.6*}
\left\{
\begin{array}{rl}
A[E_{22}-E_{11}](y_1,\dots,y_N,x_1,\dots,x_k) =& BE_{21}(y_1,\dots,y_N,x_1,\dots,x_k)\\
 -AE_{12}(y_1,\dots,y_N,x_1,\dots,x_k)  =&    CE_{21}(y_1,\dots,y_N,x_1,\dots,x_k)\\
-BE_{12}(y_1,\dots,y_N,x_1,\dots,x_k) =& C[E_{22}-E_{11}](y_1,\dots,y_N,x_1,\dots,x_k).\nonumber
\end{array}
\right.
\end{equation}
\end{definition}

\bigskip
It is verified directly  from remark \ref{rmk2.1} and the equality $E_{11}E_{22}-E_{12}E_{21}=(-1)^k$,
that 
\begin{equation}\label{eq2.6}
CE_{21}^2-BE_{21}E_{22}+AE_{22}^2=(-1)^kA. 
\end{equation}
\begin{definition}
By the Fermat-Pell conic $\mathscr{P}$ one understands the plane curve:
\begin{equation}\label{eq2.7}
Cu^2-Buv+Av^2=(-1)^kA. 
\end{equation}
\end{definition}

\begin{theorem}\label{thm2.4}
{\bf (Brock, Elkies \& Jordan \cite{BrElJo1})}
The affine variety  $V_{N,k}(\mathbf{C})$ fibers over the Fermat-Pell conic $\mathscr{P}$,
i.e. there exists a map  $\pi: V_{N,k}(\mathbf{C})\to \mathscr{P}$,   such that 
\begin{equation}\label{eq2.8}
\pi(y_1,\dots,y_N,x_1,\dots,x_k)=(E_{21}, E_{22}). 
\end{equation}
\end{theorem}

\subsection{Elliptic surfaces}

\subsubsection{Surfaces}
An algebraic surface $S$ is a variety of dimension two. 
An   elliptic surface  $S$ over a curve $C$ is a smooth projective surface
with an elliptic fibration over $C$, i.e. a surjective morphism 
$f:S\to C$  such that  almost all fibers are smooth elliptic curves.

\subsubsection{Blow-ups}
The map $\phi: S\dashrightarrow S'$ is called rational,  if it is given  
by a rational function defined everywhere except for the poles of $\phi$. 
The map $\phi$ is birational,  if the inverse $\phi^{-1}$ is a rational map. 

A birational map  $\epsilon: S\dashrightarrow S'$ is called a blow-up,
if it is defined everywhere except for a point $p\in S$ and a rational curve $C\subset S'$,
such  that  $\epsilon^{-1}(C)=p$.   
Every birational map $\phi: S\dashrightarrow S'$ is composition of a finite 
number of the blow-ups, i.e. $\phi=\epsilon_1\circ\dots\circ\epsilon_k$. 

\subsubsection{Minimal models}
The surface $S$ is called a minimal model, if any birational map $S\dashrightarrow S'$
is an isomorphism. The minimal models exist and are unique unless $S$ is 
a ruled surface. By the Castelnuovo Theorem, the surface $S$ is a minimal 
model if and only if $S$ does not contain rational curves $C$ with the 
self-intersection index $-1$.

\subsection{Non-commutative tori}
The $C^*$-algebra is an algebra  $\mathscr{A}$ over $\mathbf{C}$ with a norm 
$a\mapsto ||a||$ and an involution $\{a\mapsto a^* ~|~ a\in \mathscr{A}\}$  such that $\mathscr{A}$ is
complete with  respect to the norm, and such that $||ab||\le ||a||~||b||$ and $||a^*a||=||a||^2$ for every  $a,b\in \mathscr{A}$.  
Each commutative $C^*$-algebra is  isomorphic
to the algebra $C_0(X)$ of continuous complex-valued
functions on some locally compact Hausdorff space $X$. 
Any other  algebra $\mathscr{A}$ can be thought of as  a noncommutative  
topological space. 

The non-commutative torus  $\mathscr{A}_{\theta}$ 
is defined as a  $C^*$-algebra   generated by the unitary 
operators $u$ and $v$ satisfying  the relation $vu=e^{2\pi i\theta}uv$.
The $\mathscr{A}_{\theta}$ is said to have real multiplication (RM), if $\theta$ 
is a quadratic irrationality represented by 
the $k$-periodic continued 
fraction $[b_1,\dots, b_N;  \overline{a_1,\dots,a_k}]$. 

The non-commutative tori provide a  bridge between the 
elliptic surfaces and continued  fractions. Namely, 
there exists a  covariant functor $F$ mapping  the fibers $\mathscr{E}_t$
of an elliptic surface $S$  into the non-commutative tori   $\mathscr{A}_{\theta}$,
 see \cite[Section 1.3]{N} for the details.
The fibers  $\mathscr{E}_t$ defined over $\mathbf{Q}$ or a finite extension of $\mathbf{Q}$
map  to the non-commutative tori $F(\mathscr{E}_t)$ with RM. 

\section{Proofs}
\subsection{Proof of theorem \ref{thm1.1} }
For the sake of clarity, let us outline the main ideas.  Let $(V_{N,k}, \mathscr{P}, \pi)$ be a fiber bundle
defined by the map (\ref{eq2.8}).  
Using an exclusion process described by formulas (\ref{eq3.2}) -(\ref{eq3.4}) in below, 
we construct a  sub-bundle: 
\begin{equation}
(U_{b_1,\dots, b_N; ~a_1,\dots, a_k}, \mathscr{P}, \pi')\subset (V_{N,k}, \mathscr{P}, \pi)
\end{equation}
depending on  the point $(b_1,\dots, b_N; ~a_1,\dots, a_k)\in V_{N,k}$ and
whose  fibers are  $r$-dimensional.
Each section $\sigma$ of the bundle defines  a family of the non-commutative
tori $\mathscr{A}_{[b_1(t),\dots, b_N(t);  ~\overline{a_1(t),\dots,a_k(t)}]}$, where 
$t\in\mathscr{P}$ so that   $a_i(t_0)=a_i$ and $b_j(t_0)=b_j$. 
By the construction, 
\begin{equation}
F\left(\mathscr{E}_{\mathbf{Q}(t)}\right)=\mathscr{A}_{[b_1(t),\dots, b_N(t);  ~\overline{a_1(t),\dots,a_k(t)}]}.
\end{equation}
Let us pass to a detailed argument.

\begin{proof}
(i)   Let   $(b_1,\dots,b_N; ~a_1,\dots, a_k)$ be an integer point of 
the variety $V_{N,k}$.  To construct a subvariety  $U_{b_1,\dots,b_N; ~a_1,\dots, a_k}$, 
we denote by $(u_1,\dots, u_m)$ the variables and by $(c_1,\dots c_m)$ the constants.
Unless stated otherwise, it is assumed that 
$c_1=b_1, \dots, c_m=a_k$, where  $m=N+k$.

The equations (\ref{eq2.4})  defining the variety $V_{N,k}$ 
allow to exclude  two variables, say, $u_{m-1}$ and $u_m$,
i.e. they become algebraically dependent on the varaiables 
$\{u_i ~|~ 1\le i\le m-2\}$. Namely, one can write:
\begin{equation}\label{eq3.2}
\left\{
\begin{array}{ccc}
u_{m-1} &=& {P_{m-1}(u_1,\dots, u_{m-2})\over Q_{m-1}(u_1,\dots, u_{m-2})} \\
 u_m &=&  {P_m(u_1,\dots, u_{m-2})\over Q_m(u_1,\dots, u_{m-2})}
\end{array}
\right.
\end{equation}
for some polynomials $P_{m-1},Q_{m-1}, P_m, Q_m\in\mathbf{Z}[u_1,\dots,u_{m-2}]$. 

\medskip
Since $(c_1,\dots, c_m)\in V_{N,k}$,  one obtains  a sub-variety of the $V_{N.k}$ consisting 
of the points $(u_1,\dots, u_{m-2}, c_{m-1}, c_m)$.   In view of the (\ref{eq3.2}), 
such a sub-variety is given by the system of  equations: 
\begin{equation}\label{eq3.3}
\left\{
\begin{array}{ccc}
 P_m(u_1,\dots, u_{m-2}) &=&  c_m Q_m(u_1,\dots, u_{m-2}) \\
 P_{m-1}(u_1,\dots, u_{m-2}) &=& c_{m-1} Q_{m-1}(u_1,\dots, u_{m-2}).
\end{array}
\right.
\end{equation}

\medskip
It follows from (\ref{eq3.3}),  that again  two variables, say, $u_{m-3}$ and $u_{m-2}$
become algebraically dependent on the variables $\{u_i ~|~ 1\le i\le m-4\}$. We repeat the argument
obtaining  a sub-variety made  of the points $(u_1,\dots, u_{m-4}, c_{m-3}, c_{m-2}, c_{m-1}, c_m)$.

\medskip
It is clear, that the algorithm  will stop when the following system of the polynomial equations
in  the variables $u_i$  is  satisfied: 
\begin{equation}\label{eq3.4}
\left\{
\begin{array}{ccc}
 P_m(u_1,\dots, u_{m-2}) &=& c_m Q_m(u_1,\dots, u_{m-2}) \\
 P_{m-1}(u_1,\dots, u_{m-2}) &=&c_{m-1} Q_{m-1}(u_1,\dots, u_{m-2})\\
 &\vdots&\\
 P_2(u_1, u_2) &=& c_2 Q_2(u_1, u_2) \\
 P_1(u_1, u_2) &=& c_1 Q_1(u_1, u_2).
\end{array}
\right.
\end{equation}

\medskip
\begin{remark}
Notice that system (\ref{eq3.4}) always has a solution,  e.g. the trivial solution $(c_1,\dots, c_m)$.   
Example \ref{exm1.2}  shows that in fact such solutions can be a variety of dimension $1$. 
Below we consider a general case in terms of the Krull dimension of a polynomial ring.  
\end{remark}

\bigskip
(ii)  Let $\mathscr{I}_{c_1,\dots,c_m}$ be an ideal  generated 
by  equations (\ref{eq3.4}) in the polynomial ring $\mathbf{C}[u_1,\dots,u_m]$.  
Consider the ring  $\mathbf{C}[u_1,\dots,u_m]/\mathscr{I}_{c_1,\dots,c_m}$.
The  Krull dimension of such a ring will be denoted by $r+1$, 
see remark \ref{rmk3.2} for the notation.

\bigskip
(iii) Consider an  $(r+1)$-dimensional affine subvariety  $U_{c_1,\dots,c_m}$ of the  $V_{N,k}$
given by  equations (\ref{eq3.4}).  As it was shown in item (ii),  the points  of the variety 
$U_{c_1,\dots,c_m}$ can be written the form: 
\begin{equation}\label{eq3.5}
\left(R_1,\dots, R_{r+1},  c_{r+2},\dots, c_m\right),
\end{equation}
where $R_i\in \mathbf{Z}[u_1,\dots,u_{r+1}]$.

\bigskip
(iv)  The Brock-Elkies-Jordan Theorem \ref{thm2.4}  says that  the variety 
$V_{N,k}$ fibers over the Fermat-Pell conic $\mathscr{P}$. 
We denote by $\pi'$ a  restriction the map   $\pi: V_{N,k}\to \mathscr{P}$ to the subvariety 
 $U_{c_1,\dots,c_m}$.  Assuming  $c_1=b_1,\dots, c_m=a_k$,   one gets  a  fiber bundle 
 \begin{equation}
 (U_{b_1,\dots, b_N; ~a_1,\dots, a_k}, \mathscr{P}, \pi')
 \end{equation}
  consisting the $r$-dimensional fibers, the  $1$-dimensional base and the $(r+1)$-dimensional total space.

\bigskip
(v) Consider a global section $\sigma: \mathscr{P}\to U_{b_1,\dots, b_N; ~a_1,\dots, a_k}$. 
The $\sigma$ is given by the polynomials:  
\begin{equation}\label{eq3.6}
b_1(t), b_2(t), \dots, b_N(t)\in  \mathbf{Z}[t];  \quad a_1(t), a_2(t), \dots, a_k(t)\in \mathbf{Z}[t],
\end{equation}
 where all but $r+1$  polynomials are constants,
compare with (\ref{eq3.5}). 

\medskip
\begin{remark}\label{rmk3.2}
The number $r$ is equal to the rank of the $\mathscr{E}_{\mathbf{Q}(t)}$,
see  item (ii) of  theorem \ref{thm1.3}  to be proved in the next section.
\end{remark}

\bigskip
(vi) It remains to show,  that  
$F\left(\mathscr{E}_{\mathbf{Q}(t)}\right)=\mathscr{A}_{[b_1(t),\dots, b_N(t);  ~\overline{a_1(t),\dots,a_k(t)}]}$,
where $a_i(t)$ and $b_j(t)$ are given by formulas (\ref{eq3.6}). We shall prove this fact 
adapting  the argument of \cite{Nik1} to the case of the $\mathscr{E}_{\mathbf{Q}(t)}$. 
Namely,  suppose that the $\mathscr{E}_{\mathbf{Q}(t)}$ is given in the Legendre form:
\begin{equation}\label{eq3.7}
y^2=x(x-1)(x-\alpha(t)), \qquad\alpha(t)\in\mathbf{Q}(t). 
\end{equation}

\bigskip
(vii)  Recall that if  $F\left(\mathscr{E}_{\mathbf{Q}(t)}\right)=\mathscr{A}_{\theta}$, 
then 
\begin{equation}\label{eq3.8}
\left(\begin{matrix} b-1 & 1\cr b-2 & 1\end{matrix}\right)
\left(\begin{matrix} \theta \cr 1\end{matrix}\right)=
\left(\begin{matrix} \theta \cr 1\end{matrix}\right), 
\quad\hbox{where}\quad 
{b-2\over b+2}=\alpha(t), 
\end{equation}
see \cite[Theorem 1 \& Corollary 1.2]{Nik1}. 
Since  $b={2(1+\alpha(t))\over 1-\alpha(t)}$,  one can write (\ref{eq3.8}) in the form:
\begin{equation}\label{eq3.10}
\left(\begin{matrix} {3\alpha(t)+1\over 1-\alpha(t)} & 1\cr {4\alpha(t)\over 1-\alpha(t)} & 1\end{matrix}\right)
\left(\begin{matrix} \theta \cr 1\end{matrix}\right)=
\left(\begin{matrix} \theta \cr 1\end{matrix}\right). 
\end{equation}

\bigskip
(viii)
On the other hand, it follows from (\ref{eq2.5}) that:
\begin{equation}\label{eq3.11}
\left(\begin{matrix} E_{11} & E_{12}\cr E_{21}& E_{22}\end{matrix}\right)
\left(\begin{matrix} \theta \cr 1\end{matrix}\right)=
\left(\begin{matrix} \theta \cr 1\end{matrix}\right),
\end{equation}
 where $\theta=[b_1,\dots, b_N, \overline{a_1,\dots,a_k}]$.

\bigskip
(ix) One can factorize  matrix in (\ref{eq3.11}) as follows: 
\begin{equation}\label{eq3.12}
\begin{array}{ccc}
\left(\begin{matrix} E_{11} & E_{12}\cr E_{21}& E_{22}\end{matrix}\right)
&=&
\left(\begin{matrix} b_1 & 1 \cr 1& 0\end{matrix}\right)
\dots
\left(\begin{matrix} b_N & 1 \cr 1& 0\end{matrix}\right)
\left(\begin{matrix} a_1 & 1 \cr 1& 0\end{matrix}\right)
\dots\\
&\dots & \left(\begin{matrix} a_k & 1 \cr 1& 0\end{matrix}\right)
\left(\begin{matrix} b_N & 1 \cr 1& 0\end{matrix}\right)^{-1}
\dots
\left(\begin{matrix} b_1 & 1 \cr 1& 0\end{matrix}\right)^{-1},
\end{array}
\end{equation}
see e.g.  [Brock, Elkies \& Jordan 2021] \cite[Definition 2.4]{BrElJo1}.

 \bigskip
(x) It remains to compare  (\ref{eq3.10}) and (\ref{eq3.11}), 
i.e. 
\begin{equation}\label{eq3.13}
E_{11}=  {3\alpha(t)+1\over 1-\alpha(t)}, ~E_{21}={4\alpha(t)\over 1-\alpha(t)},
~E_{12}=E_{22}=1. 
\end{equation}
 
 \bigskip
 (xi) Since $\alpha(t)\in\mathbf{Q}(t)$, one gets $E_{12}(t), E_{21}(t)\in \mathbf{Q}(t)$. 
 Moreover, clearing the  denominators in (\ref{eq3.11}) one can always 
  assume $E_{ij}(t)\in \mathbf{Z}[t]$. 
 One obtains  from (\ref{eq3.12}),  that $a_i(t), b_j(t)\in \mathbf{Z}[t]$.

 \bigskip
 Theorem \ref{thm1.1} is proved.
 \end{proof}

\subsection{Proof of theorem \ref{thm1.3} }
\begin{proof}
Let us prove item (i) of  theorem  \ref{thm1.3}. 
 Recall that the Neron-Severi group $NS(\mathscr{E}_{\mathbf{Q}(t)})$
is the abelian group of divisors on $\mathscr{E}_{\mathbf{Q}(t)}$  modulo algebraic equivalence. 
The Picard number $\rho(\mathscr{E}_{\mathbf{Q}(t)})$ is defined as the rank of the $NS(\mathscr{E}_{\mathbf{Q}(t)})$. 
Such a number is always finite.

The idea of the proof is based on an identification of the $F(NS(\mathscr{E}_{\mathbf{Q}(t)}))$ with the 
convergents of the continued fraction  $[b_1(t),\dots, b_N(t);  ~\overline{a_1(t),\dots,a_k(t)}]$. 
By an elementary property of the continued fractions,  all such convergents are rational functions  of the first $N+k$  convergents. 
Let us pass to a detailed argument.

\bigskip
(i) Consider a global  section $\sigma_i : \mathscr{P}\to \mathscr{E}_{\mathbf{Q}(t)}$
of the elliptic surface $\mathscr{E}_{\mathbf{Q}(t)}$ with the base curve  $\mathscr{P}$. 
The $\sigma_i(\mathscr{P}):=\mathscr{P}_i$ is a genus zero curve  on  the surface  $\mathscr{E}_{\mathbf{Q}(t)}$. 
Thus one can  identify $\mathscr{P}_i$  with a divisor of the  $\mathscr{E}_{\mathbf{Q}(t)}$. 

\bigskip
(ii)  Denote by $\{\mathscr{A}_{{p\over q}} ~|~p,q\in\mathbf{Z}\}$ the  non-commutative tori
 with rational values of the parameter $\theta={p\over q}$. 
 The $\mathscr{A}_{{p\over q}}$ correspond to  the degenerate elliptic curves $\mathscr{P}$,
 i.e.  $F(\mathscr{P})=\mathscr{A}_{{p\over q}}$ \cite[Section 1.3]{N}.

\bigskip
(iii)  On the other hand, any non-commutative torus $\mathscr{A}_{\theta}$ is the 
inductive limit of an ascending sequence of the $\mathscr{A}_{{p_i\over q_i}}$,
where ${p_i\over q_i}$ are  convergents of the continued fraction of $\theta$, 
see e.g. \cite[Section 3.5]{N}. 
Consider a commutative diagram in Figure 1, where $\mathscr{P}_i$ 
is the union of all divisors of  the  $\mathscr{E}_{\mathbf{Q}(t)}$ obtained as a pull back of $F$ 
as $t$ runs through  all  admissible values.

\begin{figure}
\begin{picture}(300,110)(-40,-5)
\put(20,70){\vector(0,-1){35}}
\put(130,70){\vector(0,-1){35}}
\put(45,23){\vector(1,0){60}}
\put(45,83){\vector(1,0){60}}
\put(15,20){$\mathscr{A}_{{p_i(t)\over q_i(t)}}$}
\put(115,20){$\mathscr{A}_{[b_1(t),\dots, b_N(t);  ~\overline{a_1(t),\dots,a_k(t)}]}$}
\put(17,80){$\mathscr{P}_i$}
\put(115,80){ $\mathscr{E}_{\mathbf{Q}(t)}$}
\put(50,30){\sf embedding}
\put(50,90){\sf  embedding}
\put(0,50){$F$}
\put(110,50){$F$}
\end{picture}
\caption{}
\end{figure}

\bigskip
(iv)  The ascending sequence of the rational non-commutative  tori: 
\begin{equation}\label{eq3.14}
\mathscr{A}_{{p_1(t)\over q_1(t)}}\subset \mathscr{A}_{{p_2(t)\over q_2(t)}}
\subset\dots
\end{equation}
gives rise to an infinite inclusion sequence of the divisors
$\mathscr{P}_1\subset \mathscr{P}_2\subset\dots$
It is easy to see, that
\begin{equation}\label{eq3.15}
NS(\mathscr{E}_{\mathbf{Q}(t)})=\lim_{i\to\infty}\mathscr{P}_i. 
\end{equation}

\bigskip
(v)  Let us evaluate the number of generators of the group 
$NS(\mathscr{E}_{\mathbf{Q}(t)})$.  In view of (\ref{eq3.15}), 
this question can be reduced to the number of generators 
of the sequence (\ref{eq3.14}). 
Namely, given the periodic continued fraction 
$[b_1(t),\dots, b_N(t);  ~\overline{a_1(t),\dots,a_k(t)}]$, 
how many algebraically independent convergents 
${p_i(t)\over q_i(t)}$ are there? 

\bigskip
(vi) It is easy to see, that the total number of the independent convergents is equal 
to $N+k$. The idea is  simple: we recover the $a_i(t)$ and $b_j(t)$ 
from the first $N+k$ convergents, and then  express the remaining convergents
$\left\{{p_i(t)\over q_i(t)} ~|~i>N+k\right\}$ as the rational functions of the   $a_i(t)$ and $b_j(t)$.  

Namely,  the first convergent  ${p_1(t)\over q_1(t)}$ coincides 
with the $b_1(t)$. The second convergent ${p_2(t)\over q_2(t)}={p_1(t)\over q_1(t)}+{1\over b_2(t)}$
and, therefore, $b_2(t)=\left({p_2(t)\over q_2(t)}-{p_1(t)\over q_1(t)}\right)^{-1}$. 
Similarly, one gets $b_3(t)$ as a rational function of the ${p_1(t)\over q_1(t)}, {p_2(t)\over q_2(t)}$ and ${p_3(t)\over q_3(t)}$. 
Finally,  the $a_k(t)$ can be written as a rational function of the convergents  $\left\{{p_i(t)\over q_i(t)} ~|~1\le i\le N+k\right\}$.

Clearly,  the remaining convergents $\left\{{p_i(t)\over q_i(t)} ~|~i>N+k\right\}$ depend algebraically on the $a_i(t)$ and $b_j(t)$
and, therefore,  on the convergents  $\left\{{p_i(t)\over q_i(t)} ~|~1\le i\le N+k\right\}$.

\bigskip
(vii)  We conclude from  (i)-(vi),  that the free abelian group $NS(\mathscr{E}_{\mathbf{Q}(t)})$ 
has $N+k$ generators. In particular,   the Picard number of the surface $\mathscr{E}_{\mathbf{Q}(t)}$
 is given by the formula:
\begin{equation}
\rho(\mathscr{E}_{\mathbf{Q}(t)})=N+k.   
\end{equation}

\bigskip
Item (i) of theorem \ref{thm1.3} is proved.
 \end{proof}


\begin{proof}
Let us prove item (ii) of  theorem  \ref{thm1.3}. 
Roughly speaking,  to calculate the rank $r(\mathscr{E}_{\mathbf{Q}(t)})$  we need 
to know the  number of the convergents  $\left\{{p_i(t)\over q_i(t)} ~|~1\le i\le N+k\right\}$
independent of the parameter $t\in\mathscr{P}$. Those correspond to the ``horizontal'' 
and ``vertical'' divisors $\left\{ \mathscr{P}_i ~|~F(\mathscr{P}_i)=\mathscr{A}_{{p_i(t)\over q_i(t)}}\right\}$,
see [Sch\"utt \& Shioda 2019] \cite[Section 6.1]{SS} for the terminology.
We pass to a detailed argument. 

\bigskip
(i) Let $(U_{b_1,\dots,b_N; ~a_1,\dots, a_k}, \mathscr{P}, \pi')$
be the fiber bundle constructed in \ref{thm1.1} and let $\mathscr{F}$ be a fiber of the 
$(U_{b_1,\dots, b_N; ~a_1,\dots, a_k}, \mathscr{P}, \pi')$ having dimension $\dim\mathscr{F}=r$. 
As shown above, one can express $a_i(t)$ and $b_j(t)$ in terms of 
${p_{i}(t)\over q_{i}(t)}$ and write
points  (\ref{eq3.5}) of the variety $U_{b_1,\dots, b_N; a_1,\dots, a_k}$
in the form:
\begin{equation}\label{eq3.17}
\left(R_1,\dots, R_{r+1},  c_{r+2},\dots, c_{N+k}\right),
\end{equation}
where $R_i\in\mathbf{Z}\left[{p_1(t)\over q_1(t)},\dots, {p_{r+1}(t)\over q_{r+1}(t)}\right]$.

\bigskip
(ii)  Recall the  Tate-Shioda formula: 
\begin{equation}\label{eq3.18}
\rho(\mathscr{E}_{\mathbf{Q}(t)})=r(\mathscr{E}_{\mathbf{Q}(t)})+2+\sum_{v\in R}  (m_v-1), 
\end{equation}
where $R$ is the finite set of singular fibers of the  $\mathscr{E}_{\mathbf{Q}(t)}$  and $m_v$ is the 
number of  components of the fiber $v\in R$ [Shioda 1972] \cite[Corollary 1.5]{Shi1}.
In formula (\ref{eq3.18}) the number of horizontal divisors of  $\mathscr{E}_{\mathbf{Q}(t)}$
is equal $2$ and such of the vertical divisors is equal to $\sum_{v\in R}  (m_v-1)$
[Sch\"utt \& Shioda 2019] \cite[Section 6.1]{SS}. 

\bigskip
(iii)  Since  the horizontal and vertical divisors are not generic, they cannot depend on $t\in\mathscr{P}$. 
Thus the constants $c_i$  in  (\ref{eq3.17}) represent the horizontal and vertical divisors. 
Comparing with the Tate-Shioda formula (\ref{eq3.18}), one gets 
\begin{equation}
r(\mathscr{E}_{\mathbf{Q}(t)})=r=\dim\mathscr{F}. 
 \end{equation}

\bigskip
Item (ii) of theorem \ref{thm1.3} is proved.
\end{proof}



\begin{proof}
Let us prove item (iii) of  theorem  \ref{thm1.3}.   Roughly speaking, 
the minimal model corresponds to the surface with the least Picard
number among the surfaces in the birational equivalence class of the
$\mathscr{E}_{\mathbf{Q}(t)}$, see e.g. [Sch\"utt \& Shioda 2019] \cite[Section 4.5]{SS}.
But from \ref{thm1.3} (i) such a surface 
must minimize the sum $N+k$. We show that for the minimal model 
$N=k=1$,  so that $F(\mathscr{E}^{CM}_{p(t)})=\mathscr{A}_{[p(t); ~\overline{2p(t)}]}$
for a positive definite polynomial $p(t)\in \mathbf{Z}[t]$. 
The latter formula follows from a symmetry between the complex and 
real multiplication \cite[Section 1.4.1]{N}.   We pass to a detailed argument.

\bigskip
(i) Let  $B$ be the birational equivalence class of the   $\mathscr{E}_{\mathbf{Q}(t)}$.
Recall that a birational map $\{\mathscr{E}
\dashrightarrow \mathscr{E}' ~|~\mathscr{E}, \mathscr{E}'\in B\}$ is called dominant,
if $\rho(\mathscr{E})> \rho(\mathscr{E}')$. 
The surface $\mathscr{E}^{\min}\in B$ is a minimal model if $\mathscr{E}^{\min}$  is dominated by any other $\mathscr{E}\in B$. 

\bigskip
(ii) One gets from \ref{thm1.3} (i):  
\begin{equation}
\rho(\mathscr{E}^{\min})=\min_{\mathscr{E}\in B} ~(N+k),
\end{equation}
  where $N\ge 0$ and $k\ge 1$ satisfy (\ref{eq1.4}).  Let us find the  minimal values
  of $N$ and $k$ separately. 
  
  \bigskip
  (iii) The minimal value of $k$ is equal to $1$, since the period of continued fraction 
  (\ref{eq1.4}) cannot vanish. Thus we have $k_{\min}=1$.
  
  \bigskip
  (iv)   The minimal value of $N$ is equal to $0$, since  continued fraction   
 (\ref{eq1.4}) can be purely periodic. Thus we have $N_{\min}=0$.
 
 \bigskip
 (v) Using (iii) and (iv), one can write (\ref{eq1.4}) in the form: 
\begin{equation}\label{eq3.21}
F\left(\mathscr{E}^{\min}\right)=\left\{\mathscr{A}_{[\overline{a(t)}]}
~|~a(t)\in\mathbf{Z}[t]
\right\}.
\end{equation}

\bigskip
(vi) For an explicit construction of the  $\mathscr{E}^{\min}$, recall that
one can add a finite tail ${1\over 2}a(t)$ to the purely periodic fraction 
$[\overline{a(t)}]$. The obtained surface $(\mathscr{E}^{\min})'$ will be 
isomorphic over $\mathbf{C}$ to the original surface $\mathscr{E}^{\min}$
\cite[Theorem 1.3.1]{N}. In other words, one gets from (\ref{eq3.21}):
\begin{equation}\label{eq3.22}
F\left((\mathscr{E}^{\min})'\right)=\left\{\mathscr{A}_{[p(t); ~\overline{2p(t)}]}
~|~p(t)\in\mathbf{Z}[t]
\right\},
\end{equation}
where $p(t):={1\over 2}a(t)$. 

\bigskip
(vii)
Since $[p(t); ~\overline{2p(t)}]=(1+p^2(t))^{1\over 2}$ in (\ref{eq3.22}), one can use an explicit 
formula for the functor $F$  saying that fibers of the  surface
$(\mathscr{E}^{\min})'$ must have complex multiplication by the ring of integers of the imaginary quadratic
field  $\mathbf{Q}\left(\sqrt{-1-p^2(t)}\right)$
\cite[Theorem 1.4.1]{N}. Thus $(\mathscr{E}^{\min})'\cong \mathscr{E}^{CM}_{p(t)}$. 

\bigskip
This argument finishes the proof of item (iii) of theorem \ref{thm1.3}. 
\end{proof}


\begin{figure}
\begin{tabular}{c|c|c}
\hline
&&Picard number  \\
$D$ &  $\theta$ & of surface\\
&&$\mathscr{E}_{D(t)}$\\
\hline
$2$ &  $[1,\overline{2}]$ & $2$\\
\hline
$3$ &  $[1,\overline{1,2}]$ & $3$\\
\hline
$7$ & $[2,\overline{1,1,1,4}]$ & $5$\\
\hline
$11$ &  $[3,\overline{3,6}]$ & $3$\\
\hline
$19$ &   $[4,\overline{2,1,3,1,2,8}]$ & $7$\\
\hline
$43$ &  $[6,\overline{1,1,3,1,5,1,3,1,1,12}]$ & $11$\\
\hline
$67$ &  $[8,\overline{5,2,1,1,7,1,1,2,5,16}]$ & $11$\\
\hline
$163$ &  $[12,\overline{1, 3, 3, 2, 1, 1, 7, 1, 11, 1, 7, 1, 1, 2, 3, 3, 1, 24}]$ & $19$\\
\hline
\end{tabular}

\caption{}
\end{figure}

\section{Complex multiplication}
In this section we apply \ref{thm1.3} (i)  to calculate the Picard numbers of
 elliptic surfaces  whose fibers of rational points
(in $\mathbf{Q}P^1$)  have complex multiplication.

\bigskip
Let $D>1$ be a square-free integer. Denote by $\mathscr{E}_{CM}^{(-D)}$ 
an elliptic curve with CM by the ring of integers of the
imaginary quadratic field $\mathbf{Q}(\sqrt{-D})$. It is known, that 
\begin{equation}
F\left(\mathscr{E}_{CM}^{(-D)}\right)=\mathscr{A}_{RM}^{(D, f)},
\end{equation}
where $\mathscr{A}_{RM}^{(D, f)}$ is a non-commutative torus 
with RM  by the order of conductor $f\ge 1$ in the ring 
of integers of the real quadratic field  $\mathbf{Q}(\sqrt{D})$
\cite[Theorem 1.4.1]{N}.  In other words, one gets the $\mathscr{A}_{\theta}$, 
where 
\begin{equation}\label{eq4.2}
\theta=
\begin{cases}
\sqrt{f^2D}=[b_1; \overline{a_1, a_2,\dots, a_2, a_1, 2b_1}], & \hbox{if} ~D\equiv 2,3\mod 4\cr
 {1+\sqrt{f^2D}\over 2}=[b_1; \overline{a_1, a_2,\dots, a_2, a_1, 2b_1-1}], & \hbox{if} ~D\equiv 1\mod 4. 
\end{cases}
\end{equation}

\bigskip
As usual, denote by $(V_{1,k}(\mathbf{C}), \mathscr{P}, \pi)$ the fiber bundle corresponding to the continued
fractions (\ref{eq4.2}).  
\begin{definition}
By $\mathscr{E}_{D(t)}$ we understand an elliptic surface,  such that
for a section $\sigma: \mathscr{P}\to U_{b_1; a_1,\dots, a_k}$
of the bundle
$(U_{b_1; a_1,\dots, a_k}, \mathscr{P}, \pi')\subset (V_{1,k}(\mathbf{C}), \mathscr{P}, \pi)$
\begin{equation}
F\left(\mathscr{E}_{D(t)}\right)=
\begin{cases}
\mathscr{A}_{[b_1(t); ~\overline{a_1(t),\dots, a_1(t), 2b_1(t)}]}, & \hbox{if} ~D\equiv 2,3\mod 4\cr
\mathscr{A}_{[b_1(t); ~\overline{a_1(t),\dots, a_1(t), 2b_1(t)-1}]},  & \hbox{if} ~D\equiv 1\mod 4. 
\end{cases}
\end{equation}
\end{definition}

\bigskip
\begin{corollary}
 The Picard number of the surface $\mathscr{E}_{D(t)}$ is given by:
 \begin{equation}
 \rho(\mathscr{E}_{D(t)})=1+k,
 \end{equation}
 where $k$ is the length of period of the continued fraction (\ref{eq4.2}). 
\end{corollary}

\bigskip
\begin{example}\label{exm4.3}
Let $D=2, 3, 7, 11, 19, 43, 67$ or $163$. The imaginary quadratic field $\mathbf{Q}(\sqrt{-D})$
has class number one. Since such a  number for the real quadratic field $\mathbf{Q}(\sqrt{D})$
is also one, one gets  $f=1$  from a symmetry equation,   see  \cite[Theorem 1.4.1]{N} for the details. 
In view of $D\equiv 2, 3\mod 4$, we use the first line  in the formulas   (\ref{eq4.2}). 
 The Picard numbers of the surface $\mathscr{E}_{D(t)}$ are shown in  Figure 2. 
 Whether such surfaces are minimal is an interesting open problem.  
\end{example}
\begin{example}
Let  the surface  $\mathscr{E}_{\mathbf{Q}(t)}$ be given by equation (\ref{eq1.5}), see Example \ref{exm1.2}.  
Unlike \ref{exm4.3},  there is no complex multiplication on the fibers in this case. 
However, one gets from (\ref{eq1.6})   $N=1$ and $k=2$.  Therefore  the Picard number 
$\rho(\mathscr{E}_{\mathbf{Q}(t)})=3$. 
\end{example}

\bibliographystyle{amsplain}


\end{document}